\documentclass{amsart}
\usepackage[top=25mm, bottom=30mm, left=30mm, right=30mm]{geometry}
\usepackage[english]{babel}
\usepackage{graphicx}
\usepackage{bm}
\usepackage{amsmath,amssymb,amscd}
\usepackage{physics}
\usepackage[new]{old-arrows}
\usepackage{braket} 
\usepackage{mathrsfs}
\usepackage{tikz-cd}
\usepackage{dsfont}
\usetikzlibrary{arrows, arrows.meta, cd}
\begin{document}
\title{An Approximate-Diagonal Characterization of Amenable Actions on C*Algebras}
\author{Zhengyu Fu}
\date{}
\maketitle
\begin{abstract}
We introduce relative $h$-diagonals for crossed products $A\rtimes_{\alpha}\Gamma$
by discrete groups, which are positive-form approximate diagonals in the Haagerup
tensor product taken relative to $A$. We show that $\alpha$ is amenable if and
only if $A\rtimes_{\alpha}\Gamma$ admits a relative $h$-diagonal, an
analogue of Johnson's theorem for $L^{1}(G)$.
\end{abstract}
\begingroup
\small
\setlength{\parskip}{0pt}
\tableofcontents
\endgroup
\section{Introduction}
A Banach algebra is amenable if and only if it admits a bounded approximate diagonal \cite{Joh72}. Johnson's theorem states that $L^{1}(G)$ is amenable if and only if $G$ is amenable, for any locally compact group $G$. For C*algebras, amenability is equivalent to nuclearity \cite{Con78,Haa83}. Moreover, by Haagerup's theorem, the
diagonal may be taken of positive form $\sum_{j}d_{j}^{*}\otimes d_{j}$, which
leads to a characterization of nuclearity by positive-form approximate diagonals
in the Haagerup tensor product (Theorem 3.2, cf.~\cite{Oza02}). Recall also that the group C*algebra $C^{*}\Gamma$ of a discrete group $\Gamma$ is nuclear if and
only if $\Gamma$ is amenable.\\
~~~The aim of this paper is to establish an analogue of Johnson's theorem for C*dynamical systems
$\alpha\colon\Gamma\curvearrowright A$ of discrete groups. The group algebra is
replaced by the crossed product $A\rtimes_{\alpha}\Gamma$, and amenability of
$\Gamma$ is weakened to amenability of the action $\alpha$. The algebra side has
to be modified as well, since the naive analogue ``$A\rtimes_{\alpha}\Gamma$ admits an
approximate diagonal if and only if $\alpha$ is amenable'' fails. Indeed, for C*algebras, admitting an approximate diagonal means being nuclear, and $A\rtimes_{\alpha}\Gamma$ is
nuclear if and only if $A$ is nuclear and $\alpha$ is amenable \cite{AD87}. For instance, every action of the trivial group is amenable, whereas
$A\rtimes\{e\}=A$ admits an approximate diagonal only when $A$ is nuclear. We therefore introduce relative $h$-diagonals (Definition 3.3).\\
~~~~~~For a discrete amenable group $\Gamma$, we show that $A\rtimes_{\alpha}\Gamma$ admits
a relative $h$-diagonal (Theorem 3.4). As an application, we show that for a
$\Gamma$-C*algebra $(A,\alpha)$, if there is a separable $\alpha^{\infty}$-invariant
C*subalgebra $B\subset A^{\infty}$ that is nondegenerate on $A$ and such that
$B\rtimes_{\alpha^{\infty}|_{B}}\Gamma$ is nuclear, then $A\rtimes_{\alpha}\Gamma$ admits a
relative $h$-diagonal. Our main result is the following.\\
\\
{\bf{Theorem 1.1}}~(Theorem 5.1)\textbf{.}~For a discrete group action
$\alpha:\Gamma\curvearrowright A$, the following conditions\\
(1)~the action $\alpha$ is amenable;\\
(2)~$A\rtimes_{\alpha}\Gamma$ admits a relative $h$-diagonal\\
are equivalent.\\
\\
\section{Preliminaries}
Recall that a bounded net $(D_\lambda)_\lambda$ in the projective tensor product
$A\hat{\otimes}A$ of a Banach algebra $A$ is an approximate diagonal if
$\|a.D_\lambda-D_\lambda.a\|\to0$ and $\Delta(D_\lambda)a\to a$ for all $a\in A$,
where $a.(x\otimes y).b:=ax\otimes yb$ and $\Delta(x\otimes y):=xy$; $A$ is
amenable if and only if it admits one \cite{Joh72}, and a C*algebra is amenable
if and only if it is nuclear \cite{Con78,Haa83}.\\
\\
{\bf{The Haagerup tensor product.}}~Let $A$ be a C*algebra. For a
finite family $x=(x_{j})_{j}$ in $A$ put $\|x\|_{2}:=\|\sum_{j}x_{j}^{*}x_{j}\|^{1/2}$. Then, for finite families $x,y$ in $A$ and $c,d\in M(A)$,
\begin{equation}
  \Bigl\|\sum_{j}x_{j}^{*}y_{j}\Bigr\|\le\|x\|_{2}\,\|y\|_{2},\qquad
  \|(cx_{j}d)_{j}\|_{2}\le\|c\|\,\|x\|_{2}\,\|d\|.\nonumber
\end{equation}
We equip $A\otimes_{\mathrm{alg}}A^{\mathrm{op}}$ (with multiplication $(a\otimes b)(a'\otimes b')=aa'\otimes b'b$) with the $M(A)$-bimodule structure $c.(a\otimes b).d:=ca\otimes bd$ and the
multiplication map $\Delta(a\otimes b):=ab$; these are well defined by the
universal property of $\otimes_{\mathrm{alg}}$, since $A$ is an ideal of $M(A)$. The Haagerup norm of $w\in A\otimes_{\mathrm{alg}}A^{\mathrm{op}}$
is
$$
  \|w\|_{\mathrm h}:=\inf\Bigl\{\|(a_{j}^{*})_{j}\|_{2}\,\|(b_{j})_{j}\|_{2}:\
  w=\sum_{j}a_{j}\otimes b_{j}\Bigr\}.
$$
Let $A\otimes_{\mathrm{h}}A$ be the completion of $A\otimes_{\mathrm{alg}}A^{\mathrm{op}}$
with respect to $\|\cdot\|_{\mathrm{h}}$ (\textbf{we omit ``op'' from the notation}). From the definition and the above inequalities, one has
$$
  \|uv\|_{\mathrm{h}}\le\|u\|_{\mathrm{h}}\|v\|_{\mathrm h},\qquad
  \|c.w.d\|_{\mathrm{h}}\le\|c\|\,\|w\|_{\mathrm{h}}\,\|d\|,\qquad
  \|\Delta(w)\|\le\|w\|_{\mathrm{h}}
$$
for $u,v,w\in A\otimes_{\mathrm{alg}}A^{\mathrm{op}}$ and $c,d\in M(A)$. Hence these operations extend by continuity to $A\otimes_{\mathrm{h}}A$, which is thus a Banach
$M(A)$-bimodule, and $\Delta$ extends to a contractive bimodule map
$A\otimes_{\mathrm{h}}A\to A$.\\
For subspaces $\mathcal{S}\subset\mathcal{E}$ of a C$^*$-algebra and $w\in\mathcal{S}\otimes_{\mathrm{alg}}\mathcal{S}$, $\|w\|_{\mathrm{h},\mathcal{S}}=\|w\|_{\mathrm{h},\mathcal{E}}$, where $\|\cdot\|_{\mathrm{h},\mathcal{S}}$ is the infimum above over $a_{j},b_{j}\in\mathcal{S}$.\\
\\
{\bf{The $\Gamma$-C$^*$-correspondence $\ell^{2}(\Gamma,A)$}.} For a $\Gamma$-C*algebra $(A,\alpha)$, we regard $C_{c}(\Gamma,A)$ as an $A$-bimodule under $(a\xi)(h):=a\,\xi(h),~(\xi a)(h):=\xi(h)a$
with the $A$-valued inner product and the
induced $\Gamma$-action
$$
  \braket{\xi,\eta}:=\sum_{h\in\Gamma}\xi(h)^{*}\eta(h),\qquad
  (\tilde{\alpha}_{g}\xi)(h):=\alpha_g\bigl(\xi(g^{-1}h)\bigr)
  \qquad(g,h\in\Gamma).
$$
 Put $\|\xi\|_{2}:=\|\braket{\xi,\xi}\|^{1/2}$; this is the $2$-norm of the family $(\xi(h))_{h}$ above, so we use the same symbol. For $a\in A$, we have $\|\braket{\xi,\eta}\|\le\|\xi\|_{2}\,\|\eta\|_{2},~\|a\xi\|_{2}\le\|a\|\,\|\xi\|_{2},~\|\xi a\|_{2}\le\|\xi\|_{2}\,\|a\|$. A computation
gives
\begin{equation}
  \tilde{\alpha}_{g}(a\xi)=\alpha_{g}(a)\,(\tilde{\alpha}_{g}\xi),\qquad\tilde{\alpha}_{g}(\xi a)=(\tilde{\alpha}_{g}\xi)\,\alpha_{g}(a),\qquad
  \braket{\tilde{\alpha}_{g}\xi,\tilde{\alpha}_{g}\eta}
  =\alpha_{g}(\braket{\xi,\eta}),\nonumber
\end{equation}
in particular each $\tilde{\alpha}_{g}$ is isometric. Hence the completion
$\ell^{2}(\Gamma,A)$ of $C_{c}(\Gamma,A)$ with respect to $\|\cdot\|_{2}$, with the
continuous extensions of these operations is a
$\Gamma$-C$^*$-correspondence over $A$. For $\xi\in\ell^{2}(\Gamma,A)$ and $a\in A$ we
put $[\xi,a]:=\xi a-a\xi$.\\
\\
{\bf{Definition 2.1}}~(Amenable action)\textbf{.}~The
action $\alpha$ has the quasi-central approximation property(QAP) if
there is a net $(\xi_{i})_{i\in I}$ in $\ell^{2}(\Gamma,A)$ such that
$\|\xi_{i}\|_{2}\le1$ for all $i$, $\|[\xi_{i},a]\|_{2}\to0$ for all $a\in A$, and
$\braket{\xi_{i},\tilde{\alpha}_{g}\xi_{i}}\to1$ strictly for all $g\in\Gamma$.\\
Since $\Gamma$ is discrete, the uniform convergence on compact subsets of
$\Gamma$ required in \cite{OS21} is just pointwise convergence. By
\cite[Theorem 3.2]{OS21}, the QAP is equivalent to amenability of $\alpha$.\\
The following are well known for discrete $\Gamma$ \cite{AD87}: if $A$ is nuclear, then $\alpha$ is amenable if and only if $A\rtimes_{\mathrm{r}}\Gamma$ is nuclear; and $A\rtimes_{\mathrm{r}}\Gamma$ is nuclear if and only if $A$ is nuclear and $\alpha$ is amenable. Since nuclearity passes to quotients and $A\rtimes\Gamma=A\rtimes_{\mathrm{r}}\Gamma$ for amenable actions, the same holds for $A\rtimes\Gamma$.\\
\\
{\bf{Lemma 2.2}}~(\cite[Lemma 3.1]{CSSWW12})\textbf{.}~Let $A$ be a unital nuclear
C$^*$-algebra. For every finite set $X\Subset A$ and $\varepsilon>0$ there exist contractions
$a_{1},\dots,a_{n}\in A$ and $\theta_{1},\dots,\theta_{n}\ge0$ with $\sum_{i}\theta_{i}=1$
such that $\mu:=\sum_{i}\theta_{i}\,a_{i}^{*}\otimes a_{i}$ satisfies
$$
  \|\Delta(\mu)-1_{A}\|<\varepsilon,\qquad
  \|x.\mu-\mu.x\|_{A\hat{\otimes}A}<\varepsilon\quad(x\in X).
$$
 \section{h-Approximate Diagonal and Relative $h$-Diagonal.}
{\bf{Definition 3.1}}~(h-Approximate Diagonal)\textbf{.}~A bounded net $(D_{\lambda})_{\lambda}$ in $A\otimes_{h}A$ is an h-approximate diagonal if$$
\|a.D_{\lambda}-D_{\lambda}.a\|_{h}\to 0~,~~~~~\Delta(D_{\lambda})a\to a~~~~~~(a\in A). 
$$
A positive-form contractive h-approximate diagonal is one for which$$
D_{\lambda}=\sum_{j=1}^{m_{\lambda}}d^{*}_{\lambda,j}\otimes d_{\lambda,j}~,~~~~~~~\|\sum_{j}d^{*}_{\lambda,j}d_{\lambda,j}\|\leq 1.
$$
Since $\Delta$ is a contractive $A$-bimodule map, as in the case of an ordinary approximate diagonal, $(\Delta(D_\lambda))_\lambda$ forms a bounded approximate identity for $A$.\\
\\
{\bf{Theorem 3.2}}~(cf.\cite[Theorem 4]{Oza02})\textbf{.}~~For a C*algebra $A$, the following are equivalent.\\
(1)~$A$ is nuclear.\\
(2)~$A$ admits a positive-form contractive h-approximate diagonal.\\
\\
Proof.~For the implication (1)$\Rightarrow$(2), Haagerup's theorem gives a virtual diagonal in $\overline{\mathrm{co}}^{w*}\{x^{*}\otimes x:\|x\|\leq 1\}$; see \cite[Theorem 3.1]{Haa83}. The standard passage from a virtual diagonal to an approximate diagonal, followed by convex separation, therefore gives finite sums
$$
D_{\lambda}=\sum_{j}t_{\lambda,j}x^{*}_{\lambda,j}\otimes x_{\lambda,j},~~~~~t_{\lambda,j}\geq 0,~~\sum_{j}t_{\lambda,j}=1
,~\|x_{\lambda,j}\|\leq 1.$$
Putting $d_{\lambda,j}=t^{\frac{1}{2}}_{\lambda,j}x_{\lambda,j}$ gives the required positive-form contractive $h$-approximate diagonal, since $\|\cdot\|_{h}\leq \|\cdot\|_{\pi}$ (projective norm).\\
 For the converse, represent $A$ faithfully and nondegenerately on a Hilbert space $H$. Take $D_{\lambda}$ as in the definition, define$$
 \Phi_{\lambda}:B(H)\to B(H);~~T \mapsto \sum_{j}d^{*}_{\lambda,j}Td_{\lambda,j}.
 $$
 The maps $\Phi_{\lambda}$ are completely positive contractions. For $a\in A$,\\
 $$
 \|a\Phi_{\lambda}(T)-\Phi_{\lambda}(T)a\|\leq\|a.D_{\lambda}-D_{\lambda}.a\|_{h}\|T\|.
 $$
 Since $(\Phi_{\lambda})_{\lambda}$ are bounded linear maps from the Banach space to the von Neumann algebra, let $\Phi$ be a point-ultraweak cluster point. Then $\Phi(B(H))\subset A^{\prime}$. If $T\in A^{\prime}$, then $$
 \Phi_{\lambda}(T)=T\Delta(D_{\lambda})\to T~~~~(\mathrm{SOT}),
 $$
 and hence ultraweakly. Thus $\Phi$ is a completely positive projection of $B(H)$ onto $A^{\prime}$, which implies that $A^{\prime\prime}$ is injective. Applying this to the universal representation gives the nuclearity of $A$.~~~~$\blacksquare$\\
\\
\\
\\
{\bf{Definition 3.3}}~(Relative $h$-Diagonal)\textbf{.}~Let $A$ be a C*algebra, $C\subset A$ a nondegenerate C*subalgebra and $S\subset M(A)$ a subset such that the C*subalgebra of $M(A)$ generated by $C\cup S$ contains $A$. A net $$
 D_{\lambda}=\sum_{j=
 1}^{N_{\lambda}}d^{*}_{\lambda,j}\otimes d_{\lambda,j}\in A\otimes_{\mathrm{alg}}A
 $$
 is a relative $h$-diagonal for $(S,C)$ if, with $h_{\lambda}:=\Delta(D_{\lambda})=\sum_{j}d^{*}_{\lambda,j}d_{\lambda,j}$, the following hold.\\
 (RD1)~$(h_{\lambda})_{\lambda}$ is a contractive approximate identity of $A$.\\
 (RD2)~For every $s\in S$,  $\|s.D_{\lambda}-D_{\lambda}.s\|_{h}\to0$.\\
 (RD3)~For every $c\in C$, there are elements $c_{\lambda,1},\cdots,c_{\lambda,N_{\lambda}}\in C$ such that
 \begin{equation}
  \|(d_{\lambda,j}c-c_{\lambda,j}d_{\lambda,j})_{j}\|_{2}\to 0,~~~~~~\|(d_{\lambda,j}c^{*}-c^{*}_{\lambda,j}d_{\lambda,j})_{j}\|_{2}\to 0.\nonumber
 \end{equation}
\\
\textbf{Remark.}~There is a natural setting for crossed product C*algebras in which $A$ is a nondegenerate C*subalgebra of $A\rtimes_{\alpha}\Gamma$, the canonical unitaries $\{u_{g}\}_{g\in\Gamma}$ belong to $M(A\rtimes_{\alpha}\Gamma)$, and $\mathrm{span}_{g}\,Au_g=C_{c}(\Gamma,A)$
is a dense *subalgebra of $A\rtimes_{\alpha}\Gamma$. Clearly, for a commutative C*algebra $A$, a relative $h$-diagonal for $(A,A)$ is the same as a positive-form contractive $h$-approximate diagonal. For simplicity, we say that a crossed product $A\rtimes_\alpha\Gamma$ admits a relative $h$-diagonal (RHD) if it admits an RHD relative to $(\{u_{g}\}_{g\in\Gamma},A)$. If $A\rtimes\Gamma$ admits RHD for $(\{u_{g}\},A)$, then, passing to the quotient, $A\rtimes_{\mathrm{r}}\Gamma$ admits RHD for $(\{\lambda_{g}\},A)$.\\
\\
\\
\\
{\bf{Theorem 3.4.}}~For a C*algebra $A$ and an amenable discrete group $\Gamma$, $A\rtimes_{\alpha}\Gamma$ admits a relative $h$-diagonal.\\
\\
Proof.~Let $F_{\lambda}\Subset\Gamma$ be a F{\o}lner net and choose an approximate identity $(e_{\mu})$ of $A$ such that $\|\alpha_{g}(e_{\mu})-e_{\mu}\|\to 0$ for every $g\in \Gamma$ (\cite[Lemma 2.12]{OS21}). Put $d_{\lambda,g,\mu}:=\abs{F_{\lambda}}^{-1/2}u_{g}e_{\mu}\in A\rtimes_{\alpha}\Gamma$ and consequently $\Delta(D_{\lambda,\mu})=\abs{F_{\lambda}}^{-1}\sum_{g\in F_{\lambda}}e_{\mu}u^{*}_{g}u_{g}e_{\mu}=e_{\mu}^{2}$.\\
Since $e_{\mu}^{2}$ is also an approximate identity for $A$, (RD1) follows from the nondegeneracy of $A$.\\
(RD3) is clear, since for every $a\in A, g\in\Gamma$, $u_{g}a=\alpha_{g}(a)u_{g}.$\\
For (RD2), note that 
\begin{align*}
u_{g}.D_{\lambda,\mu}-D_{\lambda,\mu}.u_{g}&=\abs{F_{\lambda}}^{-1}\bigl(\sum_{k\in F_{\lambda}}\alpha_{g}(e_{\mu})u^{*}_{kg^{-1}}\otimes u_{k}e_{\mu}-\sum_{k\in F_{\lambda}g}e_{\mu}u^{*}_{kg^{-1}}\otimes u_{k}\alpha_{g^{-1}}(e_{\mu})\bigr)\\
&=\abs{F_{\lambda}}^{-1}\bigl(\sum_{k\in F_{\lambda}\cap F_{\lambda}g}(\alpha_{g}(e_{\mu})-e_{\mu})u^{*}_{kg^{-1}}\otimes u_{k}e_{\mu}+\sum_{k\in F_{\lambda}\cap F_{\lambda}g}e_{\mu}u^{*}_{kg^{-1}}\otimes u_{k}(e_{\mu}-\alpha_{g^{-1}}(e_{\mu}))\\
&+\sum_{k\in F_{\lambda}\backslash F_{\lambda}g}\alpha_{g}(e_{\mu})u^{*}_{kg^{-1}}\otimes u_{k}e_{\mu}-\sum_{k\in F_{\lambda}g\backslash F_{\lambda}}e_{\mu}u^{*}_{kg^{-1}}\otimes u_{k}\alpha_{g^{-1}}(e_{\mu})\bigr),
\end{align*}
gives $\|u_{g}.D_{\lambda,\mu}-D_{\lambda,\mu}.u_{g}\|_{h}\leq 2\|\alpha_{g}(e_{\mu})-e_{\mu}\|+\abs{F_{\lambda}g\triangle F_{\lambda}}/\abs{F_{\lambda}}$.~~~~$\blacksquare$\\
\\
\\
{\bf{Lemma 3.5.}}~Let $B$ be a C*algebra.\\
(1)~ If $w=(w_\lambda)$ is a net of elements of positive form with $\sup_\lambda\|w_\lambda\|_h<\infty$, then
$$
\mathcal{C}_{w}:=\{x\in M(B):\|x.w_\lambda-w_\lambda.x\|_h\to0\}
$$ is a C$^*$-subalgebra of $M(B)$.\\
(2)~Let $u=\sum_{i=1}^{M}a^{*}_{i}\otimes a_{i}\in B\otimes_{\mathrm{alg}}B$ be in positive contractive form, i.e. $h:=\Delta(u)=\sum_{i}a^{*}_{i}a_{i}$ satisfies $\|h\|\leq 1$, and let $v\in B\otimes_{\mathrm{alg}}B$ with $\|v\|_{h}\leq 1$. Then for $c,c_{1},\cdots,c_{M}\in B$, 
$$
\|c.uv-uv.c\|_{h}\leq \|(a_{i}c-c_{i}a_{i})_{i}\|_{2}+\|(a_{i}c^{*}-c^{*}_{i}a_{i})_{i}\|_{2}+M\max_{i}\|c_{i}.v-v.c_{i}\|_{h}.
$$
\\
Proof.~(1)~We only show that $\mathcal{C}_{w}$ is self-adjoint. For $v=\sum_{i}a_{i}\otimes b_{i}\in B\otimes_{\mathrm{alg}}B$, define $v^{\sharp}:=\sum_{i}b^{*}_{i}\otimes a^{*}_{i}$. For $x\in M(B)$, since $w_{\lambda}$ is of positive form, we have $(x.w_{\lambda}-w_{\lambda}.x)^{\sharp}=-(x^{*}.w_{\lambda}-w_{\lambda}.x^{*})$. Because $\|\cdot\|_{h}$ is preserved under $\sharp$, the assertion follows.
\\~~~~~~~~(2)~A direct computation gives
\begin{align*}
c.uv-uv.c&=\sum_{i}\big[(ca^{*}_{i}).v.a_{i}-a^{*}_{i}.v.(a_{i}c)\big]\\
&=\sum_{i}(ca^{*}_{i}-a^{*}_{i}c_{i}).v.a_{i}+\sum_{i}(a^{*}_{i}\otimes a_{i})(c_{i}.v-v.c_{i})+\sum_{i}a^{*}_{i}.v.(c_{i}a_{i}-a_{i}c).
\end{align*}
For the first term, regarding $(ca^{*}_{i}-a^{*}_{i}c_{i})_{i}$ as a row operator and $(a_{i})_{i}$ as a column operator, we obtain
$$
\Big\|\sum_{i}(ca^{*}_{i}-a^{*}_{i}c_{i}).v.a_{i}\Big\|_{h}\leq \|(a_{i}c^{*}-c^{*}_{i}a_{i})_{i}\|_{2}\,\|v\|_{h}\,\|(a_{i})_{i}\|_{2}.
$$
The third term is estimated similarly. Noting $\|(a_{i})_{i}\|_{2}=\|h\|^{1/2}\leq 1$ and $\|v\|_{h}\leq 1$, the desired inequality follows. ~~~~$\blacksquare$\\
\\
{\bf{Proposition 3.6.}}~Let $A,C$ and $S$ be as in Definition 3.3. If $C$ is nuclear and $A$ has a relative $h$-diagonal for $(S,C)$, then $A$ is nuclear.\\
\\
Proof.~It suffices to construct a positive-form contractive
$h$-approximate diagonal of $A$. Let $I:=\{(G,F,E,\delta):G\Subset C, F\Subset S, E\Subset A, \delta>0\}$ be directed by $(G,F,E,\delta)\prec(G',F',E',\delta')$ if
$G\subset G'$, $F\subset F'$, $E\subset E'$ and $\delta'\le\delta$.\\
We first show that for each $\iota=(G,F,E,\delta)\in I$, there is $w_{\iota}\in A\otimes_{\mathrm{alg}}A$ of contractive positive form with 
\begin{equation*}
\|\Delta(w_{\iota})x-x\|<2\delta\quad(x\in E),\qquad
\|y.w_{\iota}-w_{\iota}.y\|_h<2\delta\quad(y\in G\cup F).
\end{equation*}
Fix $\iota=(G,F,E,\delta)\in I$, choose $u=\sum_{i=1}^{M}a^{*}_{i}\otimes a_{i}\in A\otimes_{\mathrm{alg}} A$ with $\|\sum_{i}a^{*}_{i}a_{i}\|\le 1$ and
$h:=\Delta(u)$ satisfy
\[
\|hx-x\|<\delta\ \ (x\in E),\qquad \|s.u-u.s\|_h<\delta\ \ (s\in F),
\]
and for each $c\in G$ there are $c_{1}(c),\ldots,c_{M}(c)\in C$ with
\[
\|(a_{i}c-c_{i}(c)a_{i})_{i}\|_{2}+\|(a_{i}c^{*}-c_{i}(c)^{*}a_{i})_{i}\|_{2}<\delta
\]
Since $C$ is nuclear and nondegenerate in $A$, choose a positive-form contractive
h-approximate diagonal $v=\sum_{j}z^{*}_{j}\otimes z_{j}\in C\otimes_{\mathrm{alg}}C$ with $q:=\Delta(v)$,
\[
\|c_{i}(c).v-v.c_{i}(c)\|_h<\delta/M\ \ (c\in G,\ 1\le i\le M),\qquad
\Bigl\|\sum_{i}a_i^*(qa_{i}-a_{i})x\Bigr\|<\delta\ \ (x\in E).
\]
Put $w_{\iota}:=uv=\sum_{i}a^{*}_{i}.v.a_{i}=\sum_{i,j}(z_{j}a_{i})^{*}\otimes(z_{j}a_{i})$.
It is of positive form, and $\|w_{\iota}\|_{h}\le\|u\|_{h}\|v\|_{h}\le 1$. For $x\in E$
$$
\|\Delta(w_{\iota})x-x\|\leq\|\sum_{i}a^{*}_{i}(qa_{i}-a_{i})x\|+\|hx-x\| < 2\delta.
$$
Let $s\in F$, since $s.w_{\iota}=\sum_{i}(sa^{*}_{i}).v.a_{i}=(s.u)v$ and
$w_{\iota}.s=\sum_{i}a^{*}_{i}.v.(a_{i}s)=(u.s)v$, we have
$$
\|s.w_{\iota}-w_{\iota}.s\|_{h}=\|(s.u-u.s)v\|_{h}\le\|s.u-u.s\|_{h}\|v\|_{h}<\delta .
$$
Let $c\in G$, Lemma 3.5 (2) gives
$$
\|c.w_{\iota}-w_{\iota}.c\|_{h}
\le\|(a_{i}c-c_{i}(c)a_{i})_{i}\|_2+\|(a_{i}c^{*}-c_{i}(c)^{*}a_{i})_{i}\|_{2}+M\max_{i}\|c_{i}(c).v-v.c_{i}(c)\|_{h}
<2\delta .
$$
Hence $\|y.w_\iota-w_\iota.y\|_h<2\delta$ for all $y\in G\cup F$. In conclusion, we have shown that
$$
\Delta(w_{\iota})x\to x\quad(x\in A),~~~~~~
\|y.w_{\iota}-w_{\iota}.y\|_{h}\to 0\quad(y\in C\cup S).
$$
Since $\mathcal{C}_{(w_\iota)}=\{y\in M(A):\|y.w_{\iota}-w_{\iota}.y\|_{h}\to 0\}$ is a C*subalgebra of $M(A)$, and $C\cup S\subset\mathcal{C}_{(w_{\iota})}$, it contains $A$. Thus $\|x.w_{\iota}-w_{\iota}.x\|_{h}\to 0$ for all $x\in A$, so $(w_{\iota})_{\iota\in I}$ is a positive-form contractive
$h$-approximate diagonal of $A$.~~~~$\blacksquare$
\\
\\
{\bf{Corollary 3.7.}}~For a nonamenable discrete group $\Gamma$, $A\rtimes_{\alpha}\Gamma$ need not admit a relative $h$-diagonal.\\
\\
Proof.~It is well known that, for a discrete group $\Gamma$, $C^{*}\Gamma$ is nuclear if and only if $\Gamma$ is amenable. Consider $C=\mathbb{C}$ with the trivial action. If $\Gamma$ is nonamenable, then $C^{*}\Gamma$ is not nuclear, and hence, by Proposition 3.6, $\mathbb{C}\rtimes\Gamma$ does not admit a relative $h$-diagonal.~~~~$\blacksquare$\\
\\
\\
\section{Applications of RHD.}
{\bf{Lemma 4.2.}}~Let $D$ be a nuclear C$^{*}$-algebra, $D_{0}\subset D$ a dense subspace, $f\in D_{+}$ with $\|f\|\le1$, $v_{1},\ldots,v_{p}$ unitaries in $M(D)$, and $\eta>0$. Then there are $x_{1},\ldots,x_{N}\in D_{0}$ such that $m:=\sum_{i}x_{i}^{*}\otimes x_{i}$ satisfies
\begin{align*}
&\Bigl\|\sum_{i}x_{i}^{*}x_{i}\Bigr\|\le1,\qquad\Bigl\|\sum_{i}x_{i}^{*}x_{i}-f^{2}\Bigr\|<\eta,\\
&\|v_{r}.m.v_{r}^{*}-m\|_{h}<2\eta+2\|v_{r}f-fv_{r}\|\qquad(1\le r\le p).
\end{align*}
\\
Proof.~The unitization $\tilde{D}=D+\mathbb{C}1\subset M(D)$ is nuclear. By Haagerup's theorem in the form of Lemma 2.2, applied to $\tilde D$ and the finite set $\{f,v_{1}f,\ldots,v_{p}f\}\subset D$, there are contractions $\xi_{i}\in\tilde D$ and $\theta_{i}\ge0$ with $\sum_{i}\theta_{i}=1$ such that $\mu:=\sum_{i}\theta_{i}\xi_{i}^{*}\otimes\xi_{i}$ satisfies $\|\sum_{i}\theta_{i}\xi_{i}^{*}\xi_{i}-1\|<\eta$ and $\|x.\mu-\mu.x\|_{h}<\eta$ for $x\in\{f,v_{1}f,\ldots,v_{p}f\}$. Moreover $\|\mu\|_{h}\le\|\sum_{i}\theta_{i}\xi_{i}^{*}\xi_{i}\|\le1$. Put $x_{i}:=\theta_{i}^{1/2}\xi_{i}f\in D$. Then $m=f.\mu.f$ and $\sum_{i}x_{i}^{*}x_{i}=f(\sum_{i}\theta_{i}\xi_{i}^{*}\xi_{i})f$, which gives the first two estimates. For $v=v_{r}$,
\[
v.m.v^{*}-m=\bigl((vf).\mu-\mu.(vf)\bigr).(fv^{*})-(f.\mu-\mu.f).f+\mu.\bigl(vf^{2}v^{*}-f^{2}\bigr),
\]
and $\|vf^{2}v^{*}-f^{2}\|\le2\|vfv^{*}-f\|=2\|vf-fv\|$, which gives the third. By density, we can take $x_{i}$ in $D_{0}$.~~~~$\blacksquare$\\
\\
For a C*-algebra $A$, put
$$
A^{\infty}=\ell^{\infty}(\mathbb N,A)/c_{0}(\mathbb N,A),~~~~A_{\infty}=A^{\infty}\cap A^{\prime},
$$
and denote by $\alpha^{\infty}$ the induced action on $A^{\infty}$. We also say that a subalgebra $B\subset A_{\infty}$ is nondegenerate on $A\subset A^{\infty}$ in the natural way.\\
\\
{\bf{Theorem 4.4}}~\textbf{.}~Let $A$ be a separable C*algebra, let $\Gamma$ be a discrete group, and let $\alpha:\Gamma\curvearrowright A$ be an action. Suppose that there is a separable $\alpha^{\infty}$-invariant C*subalgebra $B\subset A_{\infty}$ such that\\
(1)~$B$ is nondegenerate on $A$;\\
(2)~$B\rtimes_{\alpha^{\infty}|_{B}}\Gamma$ is nuclear.\\
\\
Then $A\rtimes_{\alpha}\Gamma$ has a relative $h$-diagonal for the pair $(\{u_{g}\}_{g\in\Gamma},A)$. Consequently,\\
$$
A\rtimes_{\alpha}\Gamma~\mathrm{is~nuclear}~~~~\Longleftrightarrow~~~~A~\mathrm{is~nuclear}.
$$
\\
Proof.~Put $\beta:=\alpha^{\infty}|_{B}$ and $D:=B\rtimes_{\beta}\Gamma$, with canonical unitaries $v_{g}\in M(D)$. Denote the canonical conditional expectation by $E:D\twoheadrightarrow B$ and note that $B=E(D)$ is nuclear. Since $B$ is separable, the Choi-Effros lifting theorem gives a c.p.c.\ lift $(z_{n})_{n}\colon B\to\ell^{\infty}(\mathbb{N},A)$ of the inclusion $B\subset A^{\infty}$. Thus, for $b\in B$, $c\in A$ and $g\in\Gamma$,
$$
[(z_{n}(b))]=b,\quad[(\alpha_{g}(z_{n}(b)))]=\beta_{g}(b),\quad\|[z_{n}(b),c]\|\to 0.
$$
By Kadison-Schwarz inequality, we also have $ z_{n}(b)^{*}z_{n}(b)\le z_{n}(b^{*}b)$.\\
Write $s,t\in C_{c}(\Gamma,B)$ as $s=\sum_{k}v_{k}s_{k}$ and $t=\sum_{k}v_{k}t_{k}$,
where $s_{k}=E(v_{k}^{*}s)$ and $t_{k}=E(v_{k}^{*}t)$. Let $\Theta_{n}\colon C_{c}(\Gamma,B)\otimes_{\mathrm{alg}}C_{c}(\Gamma,B)\to(A\rtimes_{\alpha}\Gamma)\otimes_{\mathrm{alg}} (A\rtimes_{\alpha}\Gamma)$ be the linear map given by
\[
\Theta_{n}(s^{*}\otimes t)=\sum_{k}\bigl(u_{k}z_{n}(s_{k})\bigr)^{*}\otimes u_{k}z_{n}(t_{k}),
\]
which is well defined. For $T\in C_{c}(\Gamma,B)\otimes_{\mathrm{alg}}C_{c}(\Gamma,B)$ and $g\in\Gamma$, we have
$$
\|\Theta_{n}(T)\|_{h}\le\|T\|_{h}\ \text{for every }n,~~~~~~\|u_{g}.\Theta_{n}(T).u_{g}^{*}-\Theta_{n}(v_{g}.T.v_{g}^{*})\|_{h}\to 0.
$$
For the second, we may assume $T=s^{*}\otimes t$. Since $u_{k}z_{n}(b)u_{g}^{*}=u_{kg^{-1}}\alpha_{g}(z_{n}(b))$, we have
\begin{align*}
u_{g}.\Theta_{n}(T).u_{g}^{*}&=\sum_{k}\bigl(u_{kg^{-1}}\alpha_{g}(z_{n}(s_{k}))\bigr)^{*}\otimes u_{kg^{-1}}\alpha_{g}(z_{n}(t_{k})),\\
\Theta_{n}(v_{g}.T.v_{g}^{*})&=\sum_{k}\bigl(u_{kg^{-1}}z_{n}(\beta_{g}(s_{k}))\bigr)^{*}\otimes u_{kg^{-1}}z_{n}(\beta_{g}(t_{k})).
\end{align*}
Both sums are finite; moreover, $\|z_{n}\|\le1$ and
$\alpha_{g}(z_{n}(b))\approx z_{n}(\beta_{g}(b))$ for $b\in B$ give the second.\\
Now let $\Lambda$ be the set of triples $\lambda=(F,K,\varepsilon)$ with $F\Subset\mathrm{Ball}(A)$, $K\Subset\Gamma$ and $\varepsilon>0$, ordered by inclusion in $F,K$ and reverse order in $\varepsilon$. Fix $\lambda$ and put $\eta:=\varepsilon/5$. Let $f$ be a quasi-invariant contractive approximate identity of $B$ with $\|fa-a\|<\eta$ $(a\in F)$ and $\|\beta_{g}(f)-f\|<\eta$ $(g\in K)$. Since $\|v_{g}f-fv_{g}\|=\|\beta_{g}(f)-f\|$, Lemma 4.2 gives $x_{1},\ldots,x_{N}\in C_{c}(\Gamma,B)$ such that $m:=\sum_{i}x_{i}^{*}\otimes x_{i}$ satisfies
\[
\Bigl\|\sum_{i}x_{i}^{*}x_{i}\Bigr\|\le1,\qquad\Bigl\|\sum_{i}x_{i}^{*}x_{i}-f^{2}\Bigr\|<\eta,\qquad\|v_{g}.m.v_{g}^{*}-m\|_{h}<4\eta\quad(g\in K).
\]
Represent $x_{i}=\sum_{k}v_{k}x_{i,k}$ and put $d_{i,k,n}:=u_{k}z_{n}(x_{i,k})\in A\rtimes_{\alpha}\Gamma$. Then
$$
W_{n}:=\Theta_{n}(m)=\sum_{i,k}d_{i,k,n}^{*}\otimes d_{i,k,n},\qquad h_{n}:=\Delta(W_{n})=\sum_{i,k}z_{n}(x_{i,k})^{*}z_{n}(x_{i,k}).
$$
(RD1).~~Let $q:=E(\sum_{i}x_{i}^{*}x_{i})=\sum_{i,k}x_{i,k}^{*}x_{i,k}$. $0\le h_{n}\le z_{n}(q)$, hence $\|h_{n}\|\le\|q\|\le1$. Moreover $[(h_{n})]=\sum_{i,k}[(z_{n}(x_{i,k}))_{n}]^{*}\,[(z_{n}(x_{i,k}))_{n}]=q$.\\
As $E(f^{2})=f^{2}$, $\|q-f^{2}\|<\eta$; since $\|f^{2}a-a\|\le2\|fa-a\|$, for $a\in F$
$$
\limsup_{n}\|h_{n}a-a\|=\|qa-a\|\le\|q-f^{2}\|+\|f^{2}a-a\|<3\eta .
$$
(RD2).~~For $g\in K$, by $\|T.u_{g}^{*}\|_{h}=\|T\|_{h}$,
\begin{align*}
\limsup_{n}\|u_{g}.W_{n}-W_{n}.u_{g}\|_{h}&=\limsup_{n}\|u_{g}.\Theta_{n}(m).u_{g}^{*}-\Theta_{n}(m)\|_{h}\\
&\le\limsup_{n}\bigl\|\Theta_{n}(v_{g}.m.v_{g}^{*}-m)\bigr\|_{h}\le\|v_{g}.m.v_{g}^{*}-m\|_{h}<4\eta.
\end{align*}
(RD3).~~For $c\in A$, $d_{i,k,n}c-\alpha_{k}(c)d_{i,k,n}=u_{k}[z_{n}(x_{i,k}),c]\to 0$. As the index set is finite, $\|(d_{i,k,n}c-\alpha_{k}(c)d_{i,k,n})_{i,k}\|_{2}\to 0$, and likewise with $c^{*}$ in place of $c$.\\
 Since $4\eta<\varepsilon$, there is $n_{\lambda}$ such that $D_{\lambda}:=W_{n_{\lambda}}$ and $h_{\lambda}:=\Delta(D_{\lambda})$ satisfy $\|h_{\lambda}a-a\|<\varepsilon$ $(a\in F)$, $\|u_{g}.D_{\lambda}-D_{\lambda}.u_{g}\|_{h}<\varepsilon$ $(g\in K)$, and both $2$-norms in (RD3) are less than $\varepsilon$ for $c\in F$. Every element of $\mathrm{Ball}(A)$, resp.\ of $\Gamma$, eventually lies in $F$, resp.\ $K$; hence, by scaling, $h_{\lambda}a\to a$ for $a\in A$, and (RD2), (RD3) hold. Thus $(D_{\lambda})_{\lambda\in\Lambda}$ is a relative $h$-diagonal for $(\{u_{g}\}_{g\in\Gamma},A)$. If $A$ is nuclear, then so is $A\rtimes_{\alpha}\Gamma$ by Proposition 3.6, which completes the proof.~~~~$\blacksquare$\\
\\
\\
\textbf{Example.}~Let $\Gamma$ be a countable discrete group, let $\alpha\colon\Gamma\curvearrowright A$
be an action on a separable nuclear $C^*$-algebra, and let
$\mathrm{lt}\colon\Gamma\curvearrowright C_{0}(\Gamma)$ be the left translation.
It follows from Theorem 4.4 that $(A\otimes C_{0}(\Gamma))\rtimes_{\alpha\otimes\mathrm{lt}}\Gamma$ is nuclear. If $A$ is unital, we can consider $B$ as $\{[(1\otimes f)_{n}]_{n} : f\in C_{0}(\Gamma)\}\subset (A\otimes C_{0}(\Gamma))_{\infty}$ and the other conditions are clear. For nonunital $A$, by \cite[Lemma, p.~152]{Kas88} take an approximate identity $(h_{n})_{n\in\mathbb{N}}$ of $A$ such that $\|[h_{n},a]\|\to0$ and $\|\alpha_{g}(h_{n})-h_{n}\|\to0$ for all $a\in A$ and $g\in\Gamma$. Since $C_{0}(0,1]=\{f\in C[0,1]:f(0)=0\}$, we have $f(h_{n})\in A$ for $f\in C_{0}(0,1]$, so by functional calculus we can define a $*$-homomorphism
$$
  \varphi:C_{0}(0,1]\otimes C_{0}(\Gamma)\to (A\otimes C_{0}(\Gamma))_{\infty};\qquad
f\otimes x \mapsto [(f(h_{n})\otimes x)_{n}]_{n}.
$$
It is $\Gamma$-equivariant because $\|\alpha_{g}(h_{n})-h_{n}\|\to0$ and the functional calculus is continuous. Set $B:=\mathrm{Im}\,\varphi$ which is invariant, and the nondegeneracy of $B$ is routine. The equivariant surjection $\varphi$ induces a surjection
$$
    C_{0}(0,1]\otimes\mathbb{K}(\ell^2\Gamma)\simeq
    \bigl(C_0(0,1]\otimes C_0(\Gamma)\bigr)\rtimes_{\mathrm{id}\otimes\mathrm{lt}}\Gamma
    \twoheadrightarrow B\rtimes_{(\alpha\otimes\mathrm{lt})^\infty|_B}\Gamma .
$$
which implies $B\rtimes\Gamma$ is nuclear and we are done.~~~~$\blacksquare$\\
\\
\\
\textbf{Remark.}~In some cases (e.g., for amenable actions), Theorem 4.4 also tells us that for a nonamenable group $\Gamma$, $A\rtimes\Gamma$ may have relative $h$-diagonals. Let $\partial\mathbb{F}_{2}$ denote the Gromov boundary of $\mathbb{F}_{2}$ and consider the boundary action $\beta:\mathbb{F}_{2}\curvearrowright\partial\mathbb{F}_{2}$. We only need the fact that this action is amenable ($\mathbb{F}_{2}$ is exact). Consider $C(\partial\mathbb{F}_{2})$ as a constant sequences in $C(\partial\mathbb{F}_{2})^{\infty}$ and denote it by $B$. It is routine to show that $B\rtimes_{\beta^{\infty}|_{B}}\mathbb{F}_{2}
\simeq C(\partial\mathbb{F}_{2})\rtimes_{\beta}\mathbb{F}_{2}$ and hence is nuclear. Thus, $C(\partial\mathbb{F}_{2})\rtimes_{\beta}\mathbb{F}_{2}$ does have a relative $h$-diagonal.\\
\\
In view of Theorems 3.4, 4.4 and the preceding Remark, it is natural to ask whether amenability of a discrete group action can be characterized in terms of RHD. We answer this question affirmatively in the next section.\\
\\
\section{RHD and Amenability of Discrete Group Actions.}
Let $\pi_{u}:A\to B(H_{u})$ be the universal representation and put $N:=\pi_{u}(A)^{\prime\prime}\simeq A^{**}$; note that
$N^{\prime}=\pi_{u}(A)^{\prime}$. For $g\in\Gamma$ let $V_{g}$ be the unitary on $H_{u}$ given by
$V_{g}\pi_{\varphi}(a)\xi_{\varphi}:=\pi_{\varphi\circ\alpha_{g}^{-1}}(\alpha_{g}(a))\xi_{\varphi\circ\alpha_{g}^{-1}}$. Then $g\mapsto V_{g}$ is a unitary
representation with $V_{g}\pi_{u}(a)V_{g}^{*}=\pi_{u}(\alpha_{g}(a))$, so $\beta_{g}:=\mathrm{Ad}\,V_{g}$
restricts to $\alpha_{g}^{**}$ on $N$ and preserves $N^{\prime}$ and $Z(N)$; in particular
$(Z(N),\beta)\simeq (Z(A^{**}),\alpha^{**})$. For a von Neumann algebra $M$ with a $\Gamma$-action
$\beta$ we identify $\ell^{\infty}\Gamma\bar\otimes M$ with $\ell^{\infty}(\Gamma,M)$; thus
$1\otimes y$ is the constant function with value $y$, and
$((\mathrm{lt}\otimes\beta)_{g}F)(s)=\beta_{g}(F(g^{-1}s))$.\\
 On $\mathcal{K}:=\ell^{2}(\Gamma,H_{u})$ we put
$$
  (\pi(a)\xi)(s):=\pi_{u}(\alpha_s(a))\xi(s),\qquad (U_{g}\xi)(s):=\xi(sg),\qquad
  (T_{g}\xi)(s):=V_{g}\xi(g^{-1}s),
$$
and $(M_{F}\xi)(s):=F(s)\xi(s)$ for $F\in\ell^{\infty}(\Gamma,N^{\prime})$. Then $(\pi,U)$ is a
nondegenerate covariant representation; let $\sigma:A\rtimes\Gamma\to B(\mathcal{K})$ be its
integrated form, so that $\pi=\sigma|_{A}$ and $U_{g}=\bar\sigma(u_{g})$, and put
$R:=\sigma(A\rtimes\Gamma)^{\prime}=\pi(A)^{\prime}\cap\{U_{g}\}_{g}^{\prime}$. Let $\iota_{s}: H_{u}\to \mathcal{K}\simeq\ell^{2}\Gamma\otimes H_{u}$ send $\eta$ to $\delta_{s}\otimes \eta$, so that $\iota_{s}^{*}\xi=\xi(s)$.\\
\\
\noindent{\bf Remark.}~The following are checked by direct computation.
\begin{enumerate}
\item[(a)] $\iota_{e}^{*}\pi(a)=\pi_{u}(a)\iota_{e}^{*}$ and $\pi(a)\iota_{e}=\iota_{e}\pi_{u}(a)$; hence
  $\iota_e^*X\iota_e\in N'$ for $X\in\pi(A)'$.
\item[(b)] $F\mapsto M_F$ is a unital $*$-homomorphism into $\pi(A)'$, and $M_{1\otimes y}\in R$
  for $y\in N^{\prime}$.
\item[(c)] $T_{g}$ is a unitary in $R$ with $T_{g}M_{F}T_{g}^{*}=M_{(\mathrm{lt}\otimes\beta)_{g}F}$ and $T_{g}^{*}\iota_{e}=\iota_{g^{-1}}V_{g}^{*}$.
\item[(d)] $\iota_{s}=U_{s^{-1}}\iota_e$; hence $\iota_{s}^{*}X\iota_{s}=\iota_{e}^{*}X\iota_{e}$ for $X\in R$.
\end{enumerate}
~~
\\
\\
{\bf{Theorem 5.1}}(Main Theorem)\textbf{.}~Let $(A,\alpha,\Gamma)$ be a discrete C* dynamical system. The following conditions are
equivalent.\\
(1)~The action $\alpha$ has the QAP.\\
(2)~There is a $\Gamma$-conditional expectation $$\ell^{\infty}(\Gamma)\bar{\otimes} Z(A^{**})\to Z(A^{**}).$$
Here we equip $\ell^{\infty}(\Gamma)\bar{\otimes} Z(A^{**})$ with the diagonal action $\mathrm{lt}\otimes\alpha^{**}.$\\
(3)~$A\rtimes_{\alpha}\Gamma$ admits RHD.\\
\\
Since (1) and (2) are equivalent by \cite[Theorem 3.2]{OS21}, we show $(3) \Rightarrow (2)$ and $(1) \Rightarrow (3)$.\\
\\
{\bf{Lemma 5.2.}}~Let $D_{\lambda}=\sum_{j}d^{*}_{\lambda,j}\otimes d_{\lambda,j}$ be a relative
$h$-diagonal of $A\rtimes\Gamma$, put $R:=\sigma(A\rtimes\Gamma)^{\prime}$, and define
$\Phi_{\lambda}(T):=\sum_{j}\sigma(d_{\lambda,j})^{*}T\sigma(d_{\lambda,j})$ for $T\in B(\mathcal{K})$. Then each $\Phi_{\lambda}$ is a completely positive contractive $R$-bimodule map, and
\begin{enumerate}
\item[(a)] $\|[\Phi_{\lambda}(T),\pi(a)]\|\to0$ for all $T\in\pi(A)^{\prime}$ and $a\in A$.
\item[(b)] $\|[U_g,\Phi_{\lambda}(T)]\|\to0$ for all $T\in B(\mathcal{K})$ and $g\in\Gamma$.
\item[(c)] $\Phi_{\lambda}(T)\to T$ in SOT for all $T\in R$.
\end{enumerate}
Proof.~Write $k_{\lambda}=\sum_{j}d_{\lambda,j}^{*}d_{\lambda,j}$. Clearly $\Phi_{\lambda}$ is completely
positive, and $\|\Phi_{\lambda}\|=\|\Phi_{\lambda}(1)\|=\|\sigma(k_{\lambda})\|\le1$. Since elements
of $R$ commute with every $\sigma(d_{\lambda,j})$, $\Phi_{\lambda}$ is an $R$-bimodule map. \\
(a) Let $c_{\lambda,j}\in A$ be as in (RD3). Note that $\sigma(d_{\lambda,j})\pi(a)=\sigma(d_{\lambda,j}a),~\pi(a)\sigma(d_{\lambda,j})^{*}=\sigma(d_{\lambda,j}a^{*})^{*}$ and since $T$ commutes with $\pi(c_{\lambda,j})$,
\[
  [\Phi_{\lambda}(T),\pi(a)]
  =\sum_{j}\sigma(d_{\lambda,j})^{*}T\sigma(d_{\lambda,j}a-c_{\lambda,j}d_{\lambda,j})-\sum_{j}\sigma(d_{\lambda,j}a^{*}-c^{*}_{\lambda,j}d_{\lambda,j})^{*}T\sigma(d_{\lambda,j}),
\]
which gives $\|[\Phi_{\lambda}(T),\pi(a)]\|\le\|T\|\bigl(\|(d_{\lambda,j}a-c_{\lambda,j}d_{\lambda,j})_{j}\|_{2}+\|(d_{\lambda,j}a^{*}-c_{\lambda,j}^{*}d_{\lambda,j})_{j}\|_2\bigr)\to 0$.\\
(b) Consider the linear map $\theta_{T}\colon(A\rtimes\Gamma)\otimes_{\mathrm{alg}}(A\rtimes\Gamma)\to B(\mathcal{K}); x\otimes y\mapsto\sigma(x)T\sigma(y)$. Then $\|\theta_{T}(w)\|\le\|T\|\|w\|_{h}$. Since $\theta_{T}(u_{g}.D_{\lambda})=U_{g}\Phi_{\lambda}(T)$ and
$\theta_{T}(D_{\lambda}.u_{g})=\Phi_{\lambda}(T)U_{g}$, we obtain
$\|[U_{g},\Phi_{\lambda}(T)]\|\le\|T\|\,\|u_{g}.D_{\lambda}-D_{\lambda}.u_{g}\|_{h}\to 0$ by (RD2).\\
(c) For $T\in R$ we have $\Phi_{\lambda}(T)=T\sigma(k_{\lambda})$, and $\sigma(k_{\lambda})\to 1$ (SOT) by (RD1) and the nondegeneracy of $\sigma$.~~~~$\blacksquare$\\
\\
\textbf{Proof of (3) $\Rightarrow$ (2).}\\
Let $D_\lambda=\sum_{j}d^{*}_{\lambda,j}\otimes d_{\lambda,j}$ be a relative h-diagonal, and take $(\Phi_{\lambda})_{\lambda}$ as in Lemma 5.2. Take a point-weak* cluster point $\Phi$ of $\Phi_{\lambda}$(assume $w^{*}\lim\Phi_{\lambda}=\Phi$), which is completely positive. By Lemma 5.2, we know $\Psi:=\Phi|_{\pi(A)^{\prime}}: \pi(A)^{\prime}\to R$ is a ucp $R$-bimodule map with $\Psi|_{R}=\mathrm{id}$.\\
Define $E(F):=\iota_{e}^{*}\Psi(M_{F})\iota_{e}$ for $F\in\ell^{\infty}(\Gamma,N^{\prime})$. By Remark (a),(b),
$E: \ell^{\infty}(\Gamma,N^{\prime})\to N^{\prime}$ is a ucp map with $E(1\otimes y)=y$ for $y\in N^{\prime}$. By
Remark (c),(d) and the $R$-bimodularity of $\Psi$,
$$
E((\mathrm{lt}\otimes\beta)_{g}F)=\iota_e^{*}T_{g}\Psi(M_{F})T_{g}^{*}\iota_{e}
=V_{g}\,\iota_{g^{-1}}^{*}\Psi(M_{F})\iota_{g^{-1}}\,V_{g}^{*}=\beta_{g}(E(F)).
$$
Since $1\otimes N^{\prime}$ lies in the multiplicative domain of $E$, we have
$E(F(1\otimes y))=E(F)y$ and $E((1\otimes y)F)=yE(F)$ for $y\in N^{\prime}$. If
$F\in\ell^{\infty}(\Gamma,Z(N))$, then $F$ commutes with $1\otimes N^{\prime}$, so
$E(F)\in N^{\prime}\cap N^{\prime\prime}=Z(N)$. Hence the restriction of $E$ to $\ell^{\infty}(\Gamma,Z(N))$ is a
$\beta$-equivariant conditional expectation onto $Z(N)$.~~~~$\blacksquare$\\
\\
\\
Recall that if $\alpha$ has the QAP, then for finite sets $F\subset A$, $K\subset\Gamma$ and $\varepsilon>0$, there is $\xi\in C_{c}(\Gamma,A)$ with $\|\xi\|_{2}\le1$ such that $\|[\xi,a]\|_{2}<\varepsilon$ and $\|(\braket{\xi,\tilde{\alpha}_{g}\xi}-1)a\|<\varepsilon$ for $a\in F$ and $g\in K$.\\
\\
\\
\textbf{Lemma 5.3.}~Let $\zeta\in C_{c}(\Gamma,A)$, put $x_{h}:=u_{h}^{*}\zeta(h)\in A\rtimes\Gamma$
and $W_{\zeta}:=\sum_{h}x_{h}^{*}\otimes x_{h}$. Then
\begin{enumerate}
\item[(a)] $x_{h}a-\alpha_{h^{-1}}(a)x_{h}=u_{h}^{*}[\zeta,a](h)$, hence
  $\|(x_{h}a-\alpha_{h^{-1}}(a)x_{h})_{h}\|_{2}=\|[\zeta,a]\|_{2}$ for $a\in A$;
\item[(b)] $\|u_{g}.W_{\zeta}-W_{\zeta}.u_{g}\|_{h}\le 2\|\zeta\|_{2}\|\zeta-\tilde{\alpha}_{g}\zeta\|_{2}$ for
 $g\in\Gamma$.
\end{enumerate}
~\\
 Proof.~For (b), since $\zeta(h)u_{g}=u_{g}\alpha_{g^{-1}}(\zeta(h))$, note that
$$
  x_{h}^{*}\otimes x_{h}u_{g}
  =u_{g}(\tilde{\alpha}_{g^{-1}}\zeta)(g^{-1}h)^{*}u_{g^{-1}h}\otimes
  u_{g^{-1}h}^*(\tilde{\alpha}_{g^{-1}}\zeta)(g^{-1}h)\qquad(h\in\Gamma);
$$
summing over $h$ gives $W_{\zeta}.u_{g}=u_{g}.W_{\tilde\alpha_{g^{-1}}\zeta}$. Hence
$$
  u_{g}.W_{\zeta}-W_{\zeta}.u_{g}
  =u_{g}.\sum_{h}\Bigl((\zeta-\tilde{\alpha}_{g^{-1}}\zeta)(h)^{*}u_{h}\otimes x_{h}
  +(\tilde{\alpha}_{g^{-1}}\zeta)(h)^{*}u_{h}\otimes u_{h}^{*}(\zeta-\tilde{\alpha}_{g^{-1}}\zeta)(h)\Bigr).
$$
Since $u_{h}u_{h}^{*}=1$ and $\|\tilde{\alpha}_{g^{-1}}\zeta\|_{2}=\|\zeta\|_{2}$, the definition of
$\|\cdot\|_{h}$ gives the estimate. Note also that
$\|\zeta-\tilde{\alpha}_{g^{-1}}\zeta\|_{2}=\|\tilde{\alpha}_{g}\zeta-\zeta\|_{2}$, since $\tilde{\alpha}_{g}$ is
isometric.~~~~$\blacksquare$\\
\\
\\
\textbf{Proof of (1) $\Rightarrow$ (3).}\\
For $\xi\in\ell^{2}(\Gamma,A)$ and $g\in\Gamma$ we write
$\tau_{g}(\xi)=\braket{\xi,\tilde{\alpha}_{g}\xi}\in A$.\\
Let $\Lambda$ be the set of triples $\lambda=(F,K,\varepsilon)$ of finite sets $F\subset A$,
$K\subset\Gamma$ and $\varepsilon>0$ and fix $\lambda=(F,K,\varepsilon)$. There exists $f\in A_{+}$
with $\|f\|\le1$ such that $\|fa-a\|<\varepsilon$ for $a\in F\cup F^{*}$ and
$\|\alpha_{g}(f)-f\|<\varepsilon$ for $g\in K$. Applying QAP to $F\cup\{f\}$, $K\cup\{e\}$ and
$\varepsilon$, we obtain $\xi\in C_{c}(\Gamma,A)$ with $\|\xi\|_{2}\le1$ such that
\[
  \|[\xi,a]\|_{2}<\varepsilon,\quad\|(\tau_{e}-1)a\|<\varepsilon\quad(a\in F),\qquad
  \|(\tau_{g}-1)f\|<\varepsilon\quad(g\in K).
\]
Put $\zeta_{\lambda}:=\xi f$, $x_{\lambda,h}:=u_{h}^{*}\zeta_{\lambda}(h)$
$(h\in\mathrm{supp}\,\zeta_{\lambda})$, note that
$\|\zeta_{\lambda}\|_{2}\le1$. We show that $D_{\lambda}:=W_{\zeta_{\lambda}}=\sum_{h}x_{\lambda,h}^{*}\otimes x_{\lambda,h}$ is an RHD.\\
(RD1).~~$k_{\lambda}=\braket{\zeta_{\lambda},\zeta_{\lambda}}=f\tau_{e}f\le f^{2}\le1$.
For $a\in F$,
$k_\lambda a-a=f\tau_{e}(fa-a)+f(\tau_{e}a-a)+(fa-a)$, so $\|k_\lambda a-a\|<3\varepsilon$. Thus $(k_{\lambda})$ is an
approximate unit of $A\rtimes\Gamma$.\\
(RD2).~~For $g\in K$ we have
$\zeta_\lambda-\tilde\alpha_g\zeta_\lambda=(\xi-\tilde{\alpha}_{g}\xi)f
+(\tilde{\alpha}_{g}\xi)(f-\alpha_g(f))$. Since $\braket{\xi,\xi}\le1$, we have
$\braket{\xi-\tilde\alpha_{g}\xi,\xi-\tilde{\alpha}_{g}\xi}\le2-\tau_{g}-\tau_{g}^{*}$, so
$$
  \|(\xi-\tilde{\alpha}_{g}\xi)f\|_{2}^{2}\le\|f(2-\tau_{g}-\tau_{g}^{*})f\|\le2\|(\tau_{g}-1)f\|<2\varepsilon,
$$
and hence $\|\zeta_{\lambda}-\tilde{\alpha}_{g}\zeta_{\lambda}\|_{2}<\sqrt{2\varepsilon}+\varepsilon$. By
Lemma 5.3 (b), $\|u_{g}.D_{\lambda}-D_{\lambda}.u_{g}\|_{h}<2(\sqrt{2\varepsilon}+\varepsilon)$ for $g\in K$; thus $\|u_{g}.D_\lambda-D_\lambda.u_{g}\|_{h}\to0$ for all $g\in\Gamma$.\\
(RD3).~~For $a\in F$ we have $[\zeta_{\lambda},a]=\xi(fa-af)+[\xi,a]f$ and
$\|fa-af\|\le\|fa-a\|+\|fa^{*}-a^{*}\|<2\varepsilon$, so $\|[\zeta_{\lambda},a]\|_{2}<3\varepsilon$.
Given $a\in A$, put $c_{\lambda,h}:=\alpha_{h^{-1}}(a)$, so that
$c_{\lambda,h}^{*}=\alpha_{h^{-1}}(a^{*})$. By Lemma 5.3 (a), the two quantities in (RD3) are
$\|[\zeta_\lambda,a]\|_{2}$ and $\|[\zeta_{\lambda},a^{*}]\|_2$, which are less than $3\varepsilon$
whenever $a,a^{*}\in F$; hence they tend to $0$.~~~~$\blacksquare$\\
\\
\\
 \textbf{Acknowledgements}:~The author wishes to express sincere gratitude to their supervisors for their invaluable support and encouragement. Heartfelt thanks are also extended to the people from Sapporo, Hokkaido---a city dear to the author's heart---for their constant discussions and kindness.\\
\\

Graduate School of Mathematics, Kyushu University, Motooka 744,\\
Nishi-ku Fukuoka 819-0395, Japan\\
Email address: canenblluba@gmail.com
\end{document}